\documentclass[1p,sort&compress]{elsarticle}
\usepackage{hyperref}
\usepackage{lineno}
\usepackage[english]{babel}
\usepackage{textcomp}
\modulolinenumbers[5]
\usepackage{amssymb}
\usepackage{amsmath} 
\usepackage{framed} 
\usepackage{cancel}

\usepackage{booktabs} 
\usepackage{makecell}

\usepackage{multicol} 
\usepackage{nomencl} 
\makenomenclature
\renewcommand*\nompreamble{\begin{multicols}{2}}
\renewcommand*\nompostamble{\end{multicols}}
\usepackage{etoolbox}

\renewcommand\nomgroup[1]{%
  \item[\itshape 
  \ifstrequal{#1}{A}{State space}{%
  \ifstrequal{#1}{B}{Subscripts}{
  \ifstrequal{#1}{C}{Energy storage}{} }}%
]\vspace{10pt}} 

\makeatletter
\providecommand{\doi}[1]{%
  \begingroup
    \let\bibinfo\@secondoftwo
    \urlstyle{rm}%
    \href{http://dx.doi.org/#1}{%
      doi:\discretionary{}{}{}%
      \nolinkurl{#1}%
    }%
  \endgroup
}
\makeatother

\addto\captionsenglish{}

\journal{Applied Energy}

\begin{document}

\begin{frontmatter}

\title{Time series aggregation for energy system design: Modeling seasonal storage}

\author[addFZJ]{Leander Kotzur\corref{corLK}}
\author[addFZJ]{Peter Markewitz}
\author[addFZJ]{Martin Robinius}
\author[addFZJ,addCFC]{Detlef Stolten}

\cortext[corLK]{Corresponding author. Email: l.kotzur@fz-juelich.de}
\address[addFZJ]{ Institute of Electrochemical Process Engineering (IEK-3), Forschungszentrum J\"ulich GmbH, Wilhelm-Johnen-Str., 52428 J\"ulich, Germany}
\address[addCFC]{Chair for Fuel Cells, RWTH Aachen University, c/o Institute of Electrochemical Process Engineering (IEK-3), Forschungszentrum J\"ulich GmbH, Wilhelm-Johnen-Str., 52428 J\"ulich, Germany}

\begin{abstract}
The optimization-based design of renewable energy systems is a computationally demanding task because of the high temporal fluctuation of supply and demand time series. In order to reduce these time series, the aggregation of typical operation periods has become common. The problem with this method is that these aggregated typical periods are modeled independently and cannot exchange energy. Therefore, seasonal storage cannot be adequately taken into account, although this will be  necessary for energy systems with a high share of renewable generation.

To address this issue, this paper proposes a novel mathematical description for storage inventories based on the superposition of inter-period and intra-period states. Inter-period states connect the typical periods and are able to account their sequence. The approach has been adopted for different energy system configurations. The results show that a significant reduction in the computational load can be achieved also for long term storage-based energy system models in comparison to optimization models based on the full annual time series.

\end{abstract}


\end{frontmatter}

\begin{table*}[!t]   
\begin{framed}
\printnomenclature
\end{framed}
\end{table*}

%
\section{Introduction: Time series aggregation for renewable energy systems}
\label{sec:Introduction}
Designing energy systems with minimal ecologic and economic impact is a highly complex task: energy supply and demand must be balanced in time, in space, and in energy form, and the increasing number of generation, storage, and load management options leads to extremely large solution spaces where identifying optimality in technology options, placement, sizing, and operation can be daunting. Solving such problems analytically may not be feasible, instead requiring the use of mathematical programs to identify the optimal solution \cite{Banos2011}. 

\subsection{Motivation to aggregate time series}
Although  Moore's Law held for the most of the last few decades \cite{Schaller1997}, the computational tractability of these mathematical programs remains substantially limited \cite{Leyffer2016}. The size of the input data directly influences the size of the related optimization problem, and with it the requirement for processing resources. The integration of renewable energy expands this challenge because the proper modelling of these technologies is only possible with increased resolution of the temporal framework \cite{Poncelet2014,Stenzel2016,Pfenninger2017}.

Therefore, it has become necessary to systematically simplify the design problem in advance. This can be done through the aggregation of the input time series to typical operational periods. This is popular because most of the considered time series have patterns to their hourly, daily and seasonal variations. Therefore, it is reasonable to reduce redundant data until the minimal required representative data set for the problem is reached. \citet{Lythcke2016} refer to these typical periods as characteristic operation patterns. 

Different methods for the aggregation of these patterns have been proposed: For example, creating typical days by averaging time series over a day defined by month or weekday has been popular \cite{Mavrotas2008,Mehleri2013,Casisi2009,Lozano2009}. Nevertheless, this approach can lead to deviations in the results of the related optimization problem due to smoothing effects in the shape of the profiles \cite{Kotzur2018,Bahl2017,Schuetz2016}. Furthermore, individual optimization methods for the aggregation of typical periods \cite{Golling2012,Poncelet2016} or graphical methods \citep{Ortiga2011} have also been introduced. In the recent literature, cluster methods have attracted growing interest for their potential to reduce sets of time series data to a few representative periods or time steps: 
The \textit{k-mean} clustering algorithm \cite{Jain2010} is probably the most popular means of aggregating the typical periods \cite{Green2014,Bahl2016,Fazlollahi2014a, Bahl2017, Adhau2014,Elsido2017,Gabrielli2017}. Alternatively, \textit{k-medoid} clustering is either used by stating a Mixed Integer Linear Problem (MILP), which is deterministically solved to an exact solution \cite{Dominguez2011,Schuetz2016}, or by applying greedy algorithms \cite{Rager2015,Stadler2016}. Another option is the \textit{hierarchical} clustering which can be used to determine groups of candidate periods by some similarity criteria \cite{Nahmmacher2016,Merrick2016}. Nevertheless, in this case an additional method must be chosen afterwards so as to decide how the cluster is represented, e.g. its medoid.

The aggregated typical periods are then integrated into the energy system model as follows: Each period defines a closed operation time frame. The economical or ecological impact of this period is represented by magnifying by the number of times it appears in the original time series. For clustering based time series aggregation it would be the cardinal number of the cluster the period represents. The sequence of its appearance in the original time series is then disregarded. 

\subsection{Typical periods and storage modeling}

This approach is challenging because its suitability is highly specific to the considered category of energy systems. For conventional system design, it could be sufficient to reduce the dataset to a few independent time slices \cite{Bahl2017,Merrick2016}, while for a storage-based system design, at least typical days are required to incorporate intra-day storage \cite{Fazlollahi2014} or typical weeks for inter-day storage \cite{Harb2015,Nahmmacher2016}. The storage inventory is thereby limited within each typical period by a so called cyclic condition \cite{Fazlollahi2014,Harb2015,Nahmmacher2016,Renaldi2017}. This defines the storage inventory at the beginning of the typical period to be equal to the storage inventory at the end of the typical period.

Going one step further, 100\% renewable energy system designs based on fluctuating renewable energy resources, like wind and photovoltaics, require adequate seasonal storage solutions \cite{Krajaci2011,Petruschke2014,Palzer2014,Zerrahn2015,Samsatli2016}. Although, alternative approaches focus more on connecting regions in order to balance weather fluctuations and try to minimize the requirement for storage, storage should be still considered as a potential solution and therefore included into energy system design models. For the appropriate modeling and scaling of these seasonal storage, time series are required that cover a whole year. 

The representative periods described with this cyclic condition, on the other hand, are only independent sections that cannot exchange energy between them. We illustrate the drawback of this formulation for storage-based energy systems by using typical weeks to design an island system largely based on a renewable energy supply \cite{Kotzur2018}. This approach results in a significant deviation of the optimal scale of the long term storage if it is compared to the optimal result based on the full time series. As this problem would be expected \cite{Lythcke2016,Nahmmacher2016}, new methods are required to solve the issue.

\citet{Rager2015} try to overcome this by grouping all the days in a month and taking the medoid as a representative day for this month. This enables the modeling of a consecutive order of these twelve days, but it has the drawback that the diversity of days in a month are not represented \cite{Nahmmacher2016}.

With respect to modeling annual storage operations, \citet{Samsatli2016} also aggregate typical days and put them in an order. The aggregation is based on their appearance in the year as well; in their case one typical day for each quarter of the year. This leads to an insufficient representation of variability within a quarter. Nevertheless, the choice of the representative period is interesting: While the demand profiles are averaged, the typical wind profiles are chosen according to their highest intra-day variability in order to aim for a robust system design.

\citet{Renaldi2017} introduce multiple time grids for the operational optimization of an energy system which also relies on seasonal storage. This approach is popular for controlling process plants or electrical grids and makes use of the different time constants of different elements of the system considered. Elements with fast response times are modeled on a time grid with a high resolution in parallel to elements with higher inertia which are considered on a time grid with low resolution. This enables a reduction in the related optimization problem in comparison to considering all elements on the same time grid. Nevertheless, the majority of energy system technologies have a varying operation inside a day. A second time grid would only reduce the variables introduced due to the seasonal storages, but the majority of the technologies still must be modeled with the full time series. Therefore, the possibility to reduce the optimization problem is limited.

\citet{Gabrielli2017} propose two new comprehensible methods (M1 and M2) for modeling seasonal storage together with time series aggregation. The majority of the system equations are also modeled with typical days while the storage equations hold for the whole original time grid (M1), which is described by a sequence of typical days. In the second method (M2), additional all equations sets that are not directly related to binary or integer decision variables are considered on the full time grid. A system operation results where the storage states of two days of the year described by the same typical day are characterized by a similar variation of stored energy but a different value of stored energy at the beginning of each day. 


\subsection{Idea and structure of the paper}

Taking the state of the art into account, we combine the approach of describing the operation by a sequence of clustered typical periods by\citet{Gabrielli2017} with the idea of describing part of the system dynamics on a second time layer, similar to \citet{Renaldi2017}: The first layer, named the intra-period time layer, models the operation within a typical period. The second layer, the inter-period time layer, considers state changes between these periods. In consequence, also the state equations and variables can be considered only once inside each typical period and once for the transitional states in the sequence of typical periods. This reduces further the size of the optimization problems in comparison to a description on the full original time grid. 

This formulation enables the modeling of inter-period storage behavior with the typical period approach, which is especially valuable for systems relying on seasonal storage. Nevertheless, the following derivation makes this approach also transferable to system models in general, e.g. state-space models used for Model Predictive Controller, where repetitive operation conditions exist that can be aggregated a priori.

The overall states, primarily the states of charge of the storages, are described by the superposition of two sub-states. Therefore, the method is only suitable for linear state models, meaning models where the dynamic equations which connect two time steps are linear. This is the case for the majority of energy system design models  \cite{Mashayekh2017,Milan2012,Wu2017,Schuetz2017,Merkel2015,Lindberg2016,Lauinger2016, Wang2015,Rieder2014,Haikarainen2016,Kwon2016}. The cost function or constraints of other components can still be non-linear.

This paper is structured as follows: The novel method is mathematically introduced and explained in Section \ref{sec:Method}. First, it is derived for state-space systems in general, and afterwards for the specific application of an energy storage technology. In Section \ref{sec:Results}, the method is applied and validated by integrating it into the same three design optimization models of different energy systems as in \citet{Kotzur2018}. The impact on computational load and the accuracy gain of the new state description is then illustrated by comparing the resulting system designs and operations to those calculated with independent typical days, as well for the results calculated with the full time series.
Finally, Section \ref{sec:Conclusion} summarizes and draws the principal conclusions.

\section{Method: Inter-period state description}
\label{sec:Method}
In this chapter, the mathematical description of the states between the typical periods is derived. Therefore, we first introduce the general algebraic state space model for equally spaced discrete time steps in section \ref{sec:x_fulltime}. In section \ref{sec:x_sep}, we determine a formulation that separates the general states into a superposition of states within a period and states between the periods. After this, we show in Section \ref{seq:x_typperiod} how the states within a period can be described by a typical period, and how this affects the description of the states between these typical periods. In the last subsection \ref{sec:storage_form}, we transfer the formulation to the state of charge for energy storage.

\subsection{General discrete time-variant state space formulation}
\label{sec:x_fulltime}
In general, the states of a system can be defined by the vector $x_{t}$ for each discrete time step $t$. The states of the next time step $x_{t+1}$ are defined according to the explicit discrete time-variant state space formulation as 
\begin{equation}
\label{eq:x_t_1}
x_{t+1} = A_t  x_{t} +  B_t u_{t} \quad \forall  \quad  t 
\end{equation}
where $A_t$ is the system matrix, $u_{t}$ describes the control vector, $B_t$ the input matrix. The description of the output vector is not introduced in this work, since it has no impact on the next time step and will not change with the time scale separation of the states.

In case of cyclic systems, the states at the beginning of the considered time frame $x_{1}$ are identical to those at the end of the time frame $x_{T+1}$ 
\begin{equation}
\label{eq:x_T}
\begin{array}{rl}
x_{T+1} = &  x_{1} \\ 
\end{array} 
\end{equation}

Furthermore, the states are often limited with lower $x_{lb}$ and upper bounds $x_{ub}$ due to technological constraints as follows:
\begin{equation}
\label{eq:x_dimensioning}
x_{lb} \leq x_{t} \leq x_{ub} \quad \forall \quad t 
\end{equation}

For a typical period, the time steps $t$ are replaced by the time steps $g$ within a single typical period $k$. The cyclic condition is then set for each period \cite{Fazlollahi2014,Harb2015,Nahmmacher2016}. The disadvantage of this formulation is, as mentioned, that the states of typical periods are then not linked in between.

\subsection{Describing the discrete states for a sequence of periods}
\label{sec:x_sep}

To overcome this issue, an alternative formulation of the equation set is proposed that integrates the states between typical periods with the help of a few auxiliary variables. 

We assume that each original time step $t \in [1,N_t]$ is represented by a time step $g\in [1,N_g]$ within a period $i\in [1,N_i]$. At this stage, the absolute number of time steps stays the same for this formulation, such that $N_t = N_g \times N_i$. 

\subsubsection{Index modification to steps and periods}
\label{seq:PeriodicDSS}
The reformulation of state equation \ref{eq:x_t_1} with a period index $i$ and an intra-period time step index $g$ results in the following equation for the states within a period:
\begin{equation}
x_{i,g+1} = A_{i,g} x_{i,g}  +  B_{i,g} u_{i,g}  \quad \forall  \quad g, i 
\label{eq:states_IxExt}
\end{equation}
and subsequently this equation corresponding to the connections between the states of the previous period with the consecutive one:
\begin{equation}
x_{i+1,1} = x_{i,N_g+1}  \quad \forall  \quad g 
\label{eq:connect_periods}
\end{equation}.

Essentially, this is only a modification of the time index without a change to the equation system.

\subsubsection{Superposition of the discrete states}

With a superposition, we divide the states into two different time layers: The inter-period states $x^{inter}_{g}$ and intra-period states $x^{intra}_{i,g}$ where the original states are represented by their sum:
\begin{equation}
x_{i,g} = x_{g}^{inter} + x_{i,g}^{intra} \quad \forall  \quad g, i 
\end{equation}
The inter-period states describe the states at the beginning of each period:
\begin{equation}
\label{eq:inter_def}
x_{i,1} = x_{i,g}^{inter}  \quad \forall  \quad i  
\end{equation}
resulting that the intra-period states are zero at the first time step:
\begin{equation}
\label{eq:first_zero}
x_{i,1}^{intra} = 0 \quad \forall  \quad i  
\end{equation}

\begin{figure}[h]
	\centering
  \includegraphics[width=0.99\columnwidth]{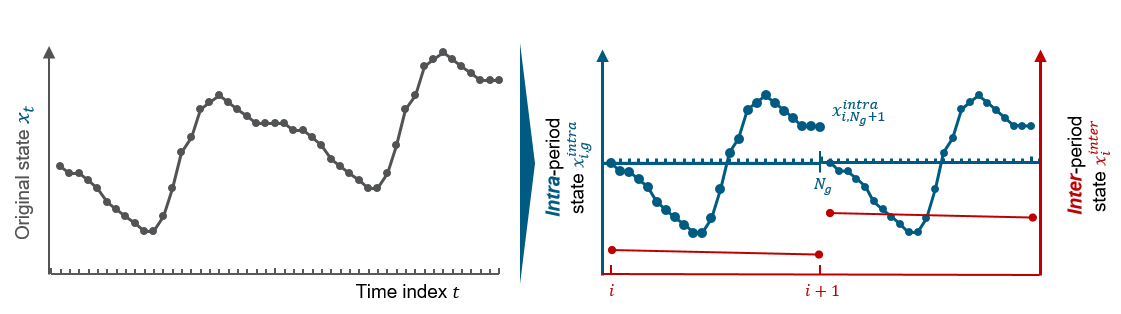}
	\caption{Separation of the original state into two states on two different time layers. These layers are here referred to as the intra- and inter period layers.}
	\label{fig:TimeScaleSeparation_sketch}
\end{figure}

The idea of this superposition is illustrated in Figure \ref{fig:TimeScaleSeparation_sketch}. 

\subsubsection{Intra-period state space equation}
Based on this superposition, we first restate the equation set from Section \ref{seq:PeriodicDSS} for the intra-period states. The equation is equivalent to \mbox{equation \ref{eq:states_IxExt}}
\begin{equation}
x^{intra}_{i,g+1} = A_{i,g} x^{intra}_{i,g}  +  B_{i,g}  u_{i,g}  
 \quad  \forall  \quad g,i 
\label{eq:intra_period}
\end{equation}

If we are not interested in the shapes of the intra-period states, we can directly enumerate for the states in the last intra-period time step $x^{intra}_{i,N_g+1}$ by a series expansion of \mbox{equation \ref{eq:intra_period}}:
\begin{equation}
\label{eq:intra_change}
x^{intra}_{i,N_g+1}  =  \quad \prod_{g} A_{i,g} x^{intra}_{i,1} \quad + \sum_{g'=1}^{N_g} \prod_{g=1}^{N_g-g'} A_{i,g}   B_{i,g} u_{i,g}   \quad  \forall  \quad i  
\end{equation}
where $x_{i,1}^{intra}$ refers to zero, as derived in equation \ref{eq:first_zero}. This is needed as an auxiliary equation in order to derive a simplified version of the description of the inter-period states.

\subsubsection{Inter-period state space equation}
\label{seq:interperiod}
Based on the introduced intra-period state equation, the inter-period state equation is derived: With the statement that the consecutive time step of the last time step of a period is the first time step of the next period, defined in equation \ref{eq:connect_periods}, following equation can also be assumed by applying a series expansion:
\begin{equation}
x_{i+1,1} = x_{i,N_g+1} = \quad \prod_{g} A_{i,g} {x_{i,1}} \quad + \sum_{g'=1}^{N_g} \prod_{g=1}^{N_g-g'} A_{i,g}   B_{i,g} u_{i,g}   \quad  \forall  \quad i 
\end{equation}
which holds for states before their superposition.

If we replace the states now with the states of the time scale separation that we defined in equations \ref{eq:inter_def} and \ref{eq:intra_change}, the following equation for the inter-period states results:
\begin{equation}
x^{inter}_{i+1} = \quad \prod_{g} A_{i,g}  x^{inter}_{i} +x^{intra}_{i,N_g+1} \quad \forall \quad i
\label{eq:inter_before}
\end{equation}
The inter period states $x^{inter}_{i+1}$ now only depend on the prior inter period states $x^{inter}_{i}$ and the final value of the intra period state $x^{intra}_{i,N_g+1}$. The advantage here is that it no longer directly depends on the input-vector $u_{i,g}$. Therefore, in the next section we can reduce all intra-period states and input variables to variables described by typical periods while still keeping the information of the states between the sequence of these periods.

\subsection{Periods to typical periods}
\label{seq:x_typperiod}
In case of an assumed aggregation of the periods to typical periods, each original candidate period $i$ belongs to a group or cluster $i \in C_k$ which is represented by its typical period $k$. In reverse, for each candidate period a typical period as $k = f(i)$ can be obtained by a look-up table.

For the case of the intra-period state equation \ref{eq:intra_period}, the result is a change of the index $i$ to index $k$. By replacing the change of the intra-period state of charge $x^{intra}_{i,N_g+1} $ of period $i$ by the change $x^{intra}_{k=f(i),N_g+1}$ of its representing period $k$, equation \ref{eq:inter_before} can be rewritten as final inter-period state equation 

\begin{equation}
x^{inter}_{s,i+1} = \quad \prod_{g} A_{i,g} x^{inter}_{i}  +  x^{intra}_{k=f(i),N_g+1} \quad \forall \quad i
\end{equation}

\begin{figure}[h]
	\centering
  \includegraphics[width=0.79\columnwidth]{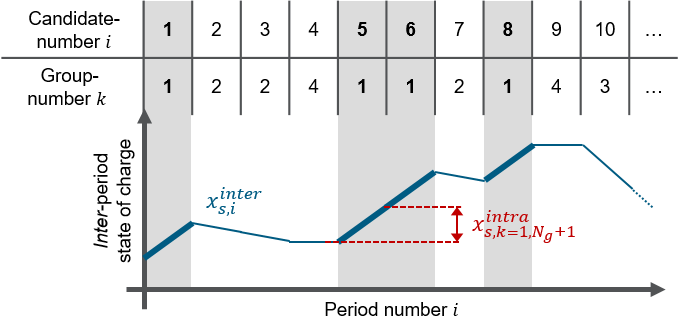}
	\caption{Sketched high layer inter-period state $x^{inter}_{i}$ based on the sequence of appearance of the representative periods $k$. This is highlighted for period or group number $1$.}
	\label{fig:InterRepresentativePeriod_sketch}
\end{figure}

The shape of the resulting inter-period states is visualized in \ref{fig:InterRepresentativePeriod_sketch}. 

Since every original state $x_{i,g}$ can be represented by the sum of the two layer states, the state constraints shown in equation \ref{eq:x_dimensioning} can be easily rewritten as 
\begin{equation}
x_{lb} \leq \quad \prod^{g}_{g'=1} A_{i,g'} x^{inter}_{i} + x^{intra}_{k=f(i),g} \leq x_{ub} \quad \forall \quad i,g
\end{equation}

The overall formulation could also be imagined in an application with more than two time layers, but this would require a more complex approach for time series aggregation and exceeds the scope of this paper.

\subsection{Storage formulation}
\label{sec:storage_form}

Since the derivation seems abstract, it is shown in the following example how these state space equations would look for the case of energy storage.

\subsubsection{Classical storage equations}
\label{sec:SOC_fulltime}

In general, the state of charge of a storage $s$ at time step $t$ can be defined by the variable $SOC_{s,t}$. With a simple Euler discretization, we can state for the state of charge in the next time step $SOC_{s,t+1}$:
\begin{equation}
\label{eq:SOC_t_1}
\begin{array}{rl}
SOC_{s,t+1} = & SOC_{s,t} (1-\eta^{self}_{s} \Delta t ) \\
& + \Delta t \left[\eta_{s}^{char} \dot{E}^{char}_{s,t} -  \frac{\dot{E}^{dis}_{s,t}}{\eta^{dis}_{s}} \right]  \\ 
\end{array} \quad \forall \quad t 
\end{equation}
where $\dot{E}^{char}_{s,t}$ describes the charging flow with an efficiency of $\eta^{char}_{s}$ and $\dot{E}^{dis}_{s,t}$ the discharging flow with related efficiency $\eta^{dis}_{s}$. $\eta^{self}_{s}$ defines the self-discharge of the storage and $\Delta t$ the step length of a single time step.

The state of charge at the beginning of the considered time frame $SOC_{s,1}$ is identical to the state of charge at the end of the time frame $SOC_{s,T+1}$  with a so-called cyclic condition:
\begin{equation}
\label{eq:SOC_T}
\begin{array}{rl}
SOC_{s,T+1} = &  SOC_{s,1} \\ 
\end{array} 
\end{equation}

The design variable of storage $s$ is described by its capacity or scale $D_s$ and limits the state of charge as follow
\begin{equation}
\label{eq:SOC_dimensioning}
0 \leq SOC_{s,t} \leq D_s \quad \forall \quad t 
\end{equation}

This set of equations (\ref{eq:SOC_t_1}, \ref{eq:SOC_T} and  \ref{eq:SOC_dimensioning}) can be found in the literature for many different storage types with slight modifications tp the syntax.

\subsubsection{Storage equations for a sequence of typical periods}
The equation set can be rewritten for the new state formulation of typical periods as follows: The intra-period states of charge are defined as 
\begin{equation}
\begin{array}{lr}
\begin{array}{ll}
SOC^{intra}_{s,k,g+1} = & SOC^{intra}_{s,k,g} (1-\eta^{self}_{s} \Delta t ) \\
& + \Delta t \left[ \eta_{s}^{char} \dot{E}^{char}_{s,k,g} -  \frac{\dot{E}^{dis}_{s,k,g}}{\eta^{dis}_{s}} \right]  
  \\
\end{array} & \quad \forall  \quad g,k  \\ 
\begin{array}{ll}
SOC^{intra}_{s,k,1}  = &0 
\end{array}
 & \quad \forall  \quad k \\
\end{array}
\end{equation}

The inter period equations are formulated as:
\begin{equation}
\begin{array}{rl}
SOC^{inter}_{s,i+1} = & SOC^{inter}_{s,i} (1-\eta^{self}_{s} \Delta t )^{N_g} + SOC^{intra}_{s,k=f(i),N_g+1} \\
 SOC^{inter}_{s,N_i+1}  = & SOC^{inter}_{s,1}
\end{array}\quad \forall \quad i 
\end{equation}

while for the limitation of the state of charge, the following equation results:
\begin{equation}
\label{eq:SOC_limit}
0 \leq  SOC^{inter}_{s,i} (1-\eta^{self}_{s} \Delta t )^{N_g}
 + SOC^{intra}_{s,k=f(i)}  \leq D_s \quad \forall  \quad i,g
\end{equation}

For small self-discharge rates, a further simplification can be used that reduces the number of equations limiting the state of charge, with the trade-off of a few additional variables. This is explained in \ref{app:soc_relaxation}.

\section{Results: Method validation and performance measure}
\label{sec:Results}
To validate the method, we apply the time series aggregation for the same three reference systems as introduced in \citet{Kotzur2018}. The overall framework, the systems and considered time series are introduced in Section \ref{sec:Framework}. In Section \ref{sec:ComparisonWithoutSequence}, we compare the optimization results of the approach with independent typical days to the newly-introduced approach with the states between the typical days and analyze the impact on the solving performance.  Finally, the storage inventory of the two-layer formulation with typical days is compared to that based on the optimization of the full time series in Section \ref{sec:StorageComp}.

\subsection{Energy system modeling framework}
\label{sec:Framework}

\begin{figure}[h]
	\centering
	 \includegraphics[width=0.99\columnwidth]{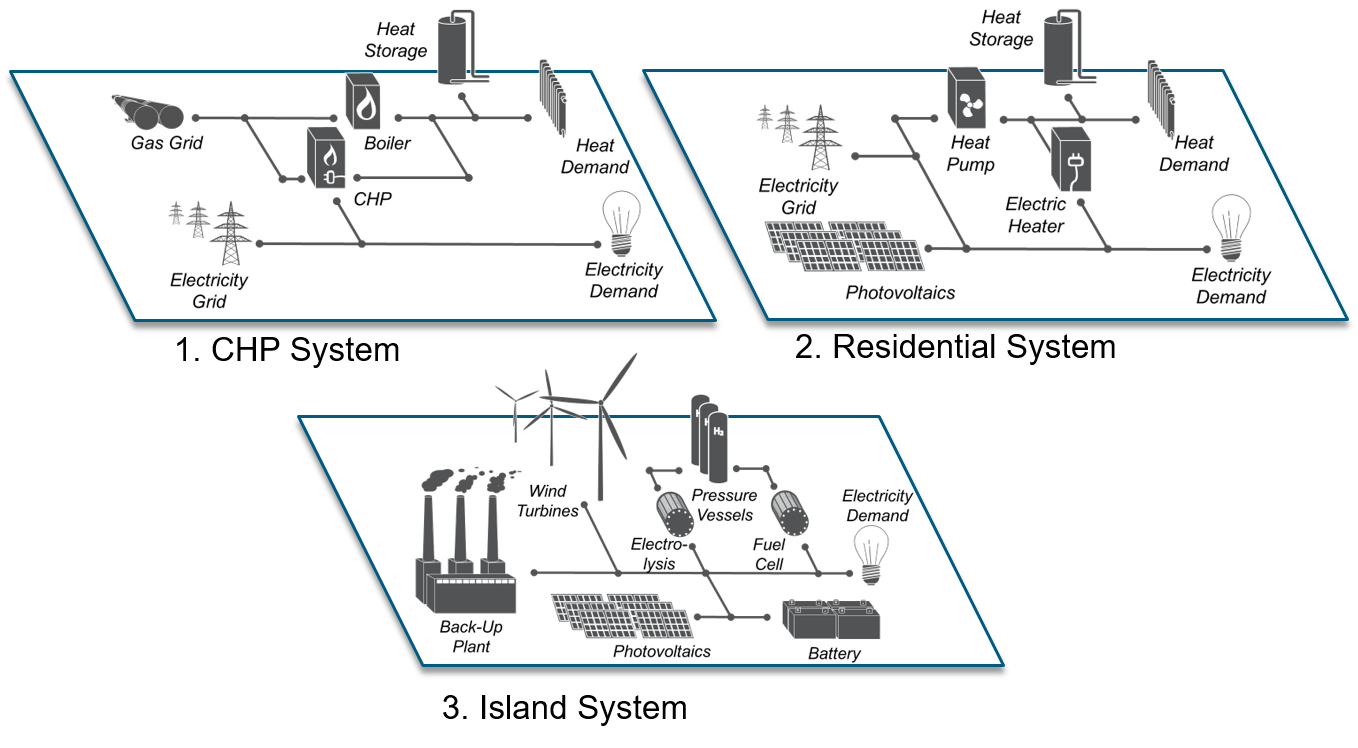}
	\caption{Structure of the three energy systems considered that are used for the validation of the method.}
	\label{fig:Systems}
\end{figure}

The newly-introduced method is applied to the following three energy supply systems:
\begin{enumerate}
\item A combined heat and power plant system (\textit{CHP}) that is supported by a peak boiler and heat storage for the supply of electricity and heat demand of a multi-family house.
The electricity demand series is gathered by down sampling the first six single residential profiles introduced by \citet{Tjaden2015} and the heat demand is simulated for a multi-family house based on a 5R1C model \cite{Schuetz2017,EN2008} with the test reference year weather data for Potsdam, Germany \cite{DWD2012} and the building data from the \textit{tabula}-database \cite{IWU2010}.
\item A \textit{residential} supply system that is based on a heatpump, an electric heater, heat storage and photovoltaics. The heat load simulation is equivalent to that for the CHP-system, but for the case of a single family house. The electricity load is the first load profile of the data from \citet{Tjaden2015}. The photovoltaic feed-in is simulated with the PV-Lib \cite{Andrews2014}.
\item An \textit{island} system that supplies the electricity for an entire region. Any transmission within the system is not considered, but the supply technologies consist of wind turbines, photovoltaics and a backup-power plant. Additionally, two storage technologies can be installed: A battery-based storage and hydrogen storage, which consists of an electrolyzer \cite{Schiebahn2015}, hydrogen pressure vessels and a fuel cell. To enforce a large share of renewable energy, the electricity supply of the backup power plant is limited to 10\% of the overall electricity consumption. The time series of the wind turbines and electricity load are drawn from \citet{Robinius2017}\cite{Robinius2015,Robinius2017a} and the photovoltaic feed-in is simulated with the PV-lib as well.
\end{enumerate}
A visualization of the final time series is seen in \ref{sec:TimeSeries}. The structure of the systems can be seen in Figure \ref{fig:Systems}.

All systems are modeled as Mixed-Integer-Linear Programs that have binary variables in the cost function of the technologies, while the operation is modeled continuous and linear. The objective function and the system constraints are found in \ref{sec:systemModel}. The techno-economic parameters are presented in \ref{sec:Parameters}.
The modeling language is \textit{Pyomo 4.3} \cite{Hart2011} and, as the solver, \textit{Gurobi 7.0.1}  \cite{Gurobi2016} was chosen. The hardware was an Intel i7-4790 CPU with 32 GB RAM, where six threads were used for the optimization.

The aggregation was performed with the \textit{tsam - Time Series Aggregation Module} \cite{Kotzur2017} where \textit{k-medoids} clustering was chosen as the aggregation method. In order to avoid side effects, no integration of the extreme periods was considered. The original time series consist of hourly data for a full year which are aggregated to typical days.

\subsection{Comparison to the typical day approach without long term states}
\label{sec:ComparisonWithoutSequence}

\begin{figure}
	\centering
  \includegraphics[width=0.99\columnwidth]{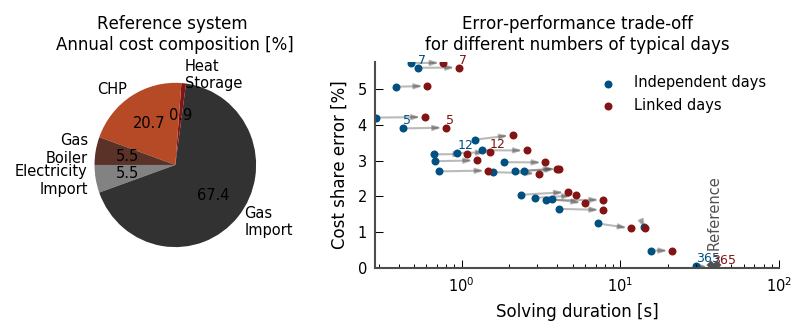}
	\caption{CHP system - Left, the cost structure of the reference optimization based on the full annual time series. Right, the cost share error for different numbers of typical days in comparison to the solving duration, once for the state formulation with independent typical periods (\textit{Independent days}), and for the state formulation that considers the sequence of the typical periods and links the states (\textit{Linked days}).}
	\label{fig:Performance_CHP}
\end{figure}

\begin{figure}[h]
	\centering
	 \includegraphics{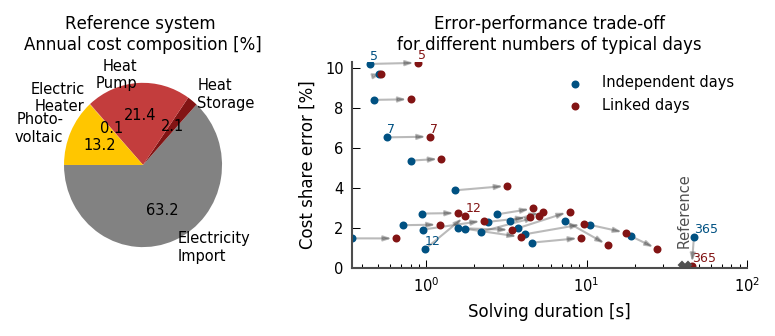}
	\caption{Residential system - Left, the cost structure of the reference optimization based on the full annual time series. Right, the cost share error for different numbers of typical days in comparison to the solving duration, once for the state formulation with independent typical periods (\textit{Independent days}), and for the state formulation that considers the sequence of the typical periods and links the states (\textit{Linked days}).}
	\label{fig:Performance_HP}
\end{figure}

\begin{figure}[h]
	\centering
 	\includegraphics[width=0.99\columnwidth]{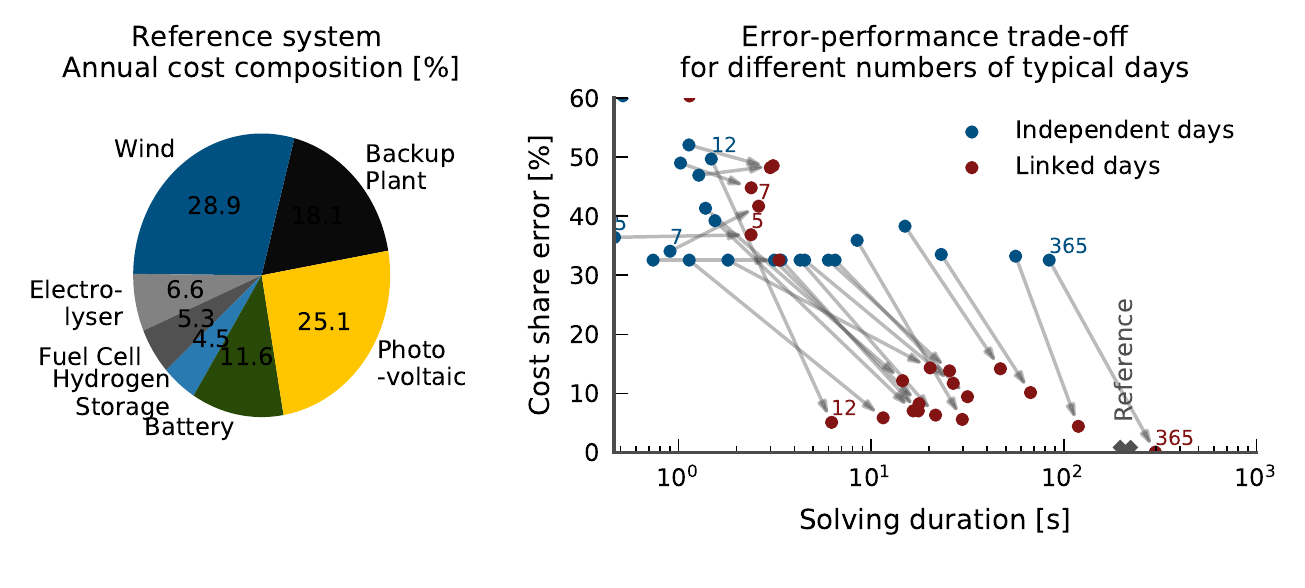}
	\caption{Island system - Left, the cost structure of the reference optimization based on the full annual time series. Right, the cost share error for different numbers of typical days in comparison to the solving duration, once for the state formulation with independent typical periods (\textit{Independent days}), and for the state formulation that considers the sequence of the typical periods and links the states (\textit{Linked days}).}
	\label{fig:Performance_Wind}
\end{figure}

The improvement or drawbacks of the new formulation are evaluated as follows: All three systems are first optimized for the full time series without aggregation. The resulting annual energy supply cost of the three systems broken down to the different technologies can be seen in the pie charts of figures \ref{fig:Performance_CHP}, \ref{fig:Performance_HP} and \ref{fig:Performance_Wind}. They determine the \textit{reference} design which should be met by the system designs based on the aggregation as accurate as possible. Then, we aggregate an increasing number of typical days using k-medoid clustering and solve the design optimization for the model formulation with independent typical days. Finally, we solve the problem with the proposed methodology where the typical days are linked in the year by the new state description. The exact design of the different systems can be found in \ref{sec:scaleResults}.

In order to have a simple performance measure to compare the results derived with the aggregated time series in comparison to the reference results, we introduce following \textit{cost share error}:
\begin{equation}
e = \frac{\sum_{d} \vert \widehat{cost}_{d} - cost_{d}  \vert }{\sum_{d} cost_{d}}
\end{equation}

This describes the deviation of the predicted device or technology cost determined with the aggregated time series $\widehat{cost}_{d}$ to the predicted cost of the reference case - the full time series - $cost_{d}$, in ratio to the full reference system cost.

The first system analyzed is the CHP system: Its cost structure can be seen in Figure \ref{fig:Performance_CHP}. The majority of the annual energy supply cost are related to energy imports in the form of electricity and gas. The integrated heat storage only counts for 0.9 \% of the overall energy cost. Therefore, the potential for improving the system results through an extension of the storage formulation is limited.

This assumption is confirmed by the results of the error measure for different typical days, seen in Figure \ref{fig:Performance_CHP} on the right. The result with independent typical days and the results with linked days, converge to an error of zero with an increasing number of typical days. The prediction error is thereby almost the same for both approaches for the same number of typical days. Nevertheless, the major difference is that the independent typical days have a faster solving performance than the approach with the linked days, while the actual difference depends on the number of typical days. This is reasonable since the inter-period state equations and variables increase the size of the optimization problem. 

For the case of the \textit{residential} system, shown in Figure \ref{fig:Performance_HP}, the heat storage has increased importance, with 2.1 \% of the overall annual energy cost. Still, this is a marginal amount compared to the other technologies and the system design does not rely on seasonal heat storage. Therefore, the comparison of the results with independent typical days and the results derived with linked days supports the conclusions made with the \textit{CHP} system: The solving performance is increased, while for a few typical days no improvement in the error measure is apparent.

The major difference is that for a high number of typical days, even up to 365 where the time series are equivalent to the original ones, the error of the approach including the link of days converges to zero, while the formulation with independent typical days retains an offset 2 \%. This highlights the major drawback of the independent typical days, which are not only a simplification of the time series, but also a restriction to the solution space of the system design since the operational possibilities are limited as well. The link of the days with the inter-period states corrects this.

The results for the last system - the \textit{island} system - can be seen in Figure \ref{fig:Performance_Wind}. Its optimal system design is primarily based on storage technologies since the major energy supply sources are wind turbines and photovoltaics. The battery accounts for 11.6 \% of the annual energy costs and the overall hydrogen storage system, including electrolyser and fuel cell, for \mbox{16.4 \%}. 

The resulting cost share errors for the system design based on the typical periods can be seen on the right of Figure \ref{fig:Performance_Wind}. First, the scale of the error is in general higher than for both systems before, since the variability of the fluctuating renewable energy is more difficult to aggregate. Therefore, more typical periods are required to achieve a system design that is similar to the system design based on the full reference data set. 

The comparison of the results with independent typical days versus those with the linked days highlights the accuracy gain of the inter-period state formulation. While for a small number of typical days both systems have a high cost share error, for more than 12 typical days the error of the formulation that includes the sequence converges to zero, although the formulation without the linkage stays at a high offset.

\begin{figure}[h]
	\centering
  \includegraphics[width=0.99\columnwidth]{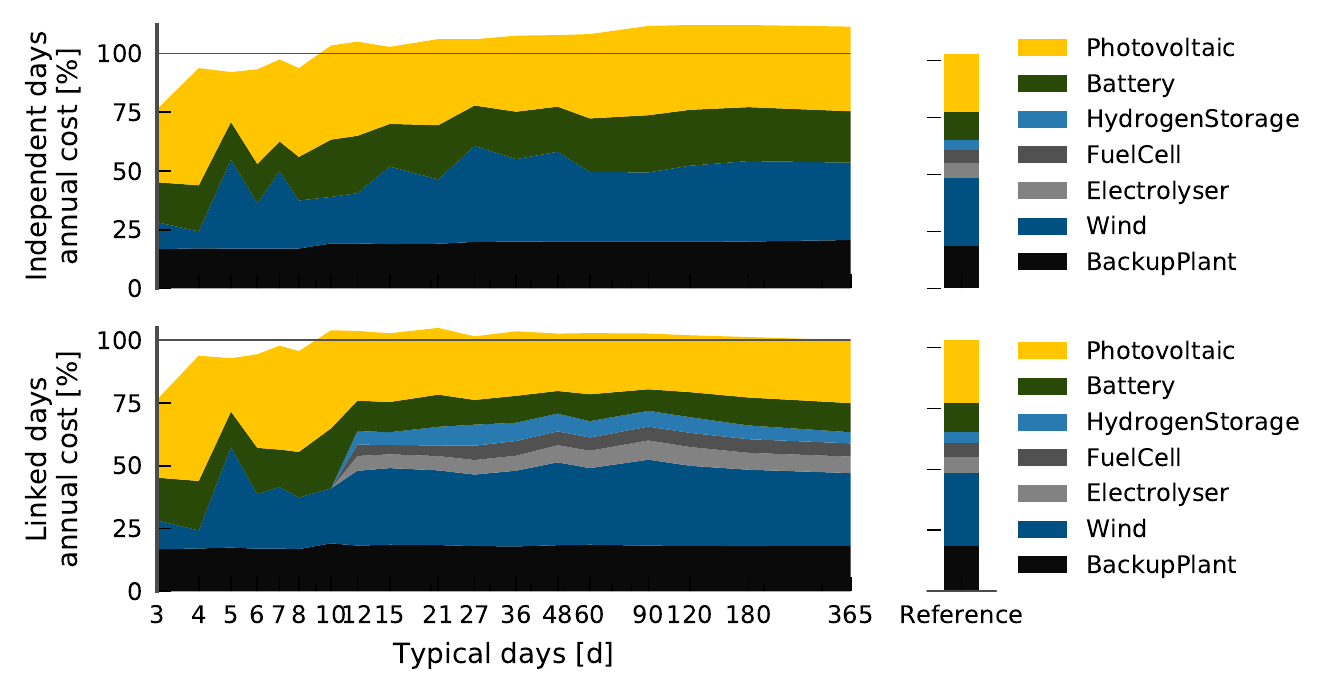}
	\caption{Comparison of the cost structure for the island system design based on independent typical days at the top and with the linked integration on the bottom; in comparison the reference system design that is based on the full time series.}
	\label{fig:CostComp}
\end{figure}

This result is shown in more detail in \mbox{Figure \ref{fig:CostComp}}. For a few typical days both approaches result in a similar system design where the battery is the only storage technology. The major reason for this is that the aggregation tends to smooth the profiles and therefore causes a reduced requirement for storage to balance the gap between renewable feed-in and electricity demand. This smoothing gets reduced with a higher number of typical periods, wherefore the amount of optimal installed storage capacities also increases. 

For the approach without the linkage of days no energy can be exchanged between the typical periods wherefore the system design must be able to supply each period independently. Therefore, also in case when many typical periods are considered, battery photovoltaic and wind must be oversized to satisfy the demand. Hydrogen storage is not cost-optimal in this case, since it is economically advantageous to the battery when energy is stored for longer operation cycles than just a single day. This drawback causes an overestimation of the necessary system cost for a high number of typical days since the operational solution space is restricted by the model.

This is different for the system formulation where the days are linked to the inter-period states: Already with 12 typical days \mbox{(6.0 sec.}, solving time) a cost structure is predicted that has a similar shape to that of the reference system \mbox{(109.2 sec.)} because the hydrogen system is included as long term storage. Nevertheless, as seen in Figure \ref{fig:CostComp}, a system design that does not significantly change for an increasing number of typical days is only found after 27 typical days  \mbox{(28.5 sec.)}. The overall annual costs are then predicted with an error of less than \mbox{2 \%} in comparison to the reference system.

This result indicates that 12 typical days could be sufficient for modeling the long term storage of hydrogen, but its accuracy still relates to the accuracy of the aggregation of typical periods itself. Therefore, the potential to improve the aggregation quality of typical periods remains. 

\subsection{Storage inventory of the island system for 12 typical days}
\label{sec:StorageComp}
\begin{figure}[h]
	\centering
  \includegraphics[width=0.99\columnwidth]{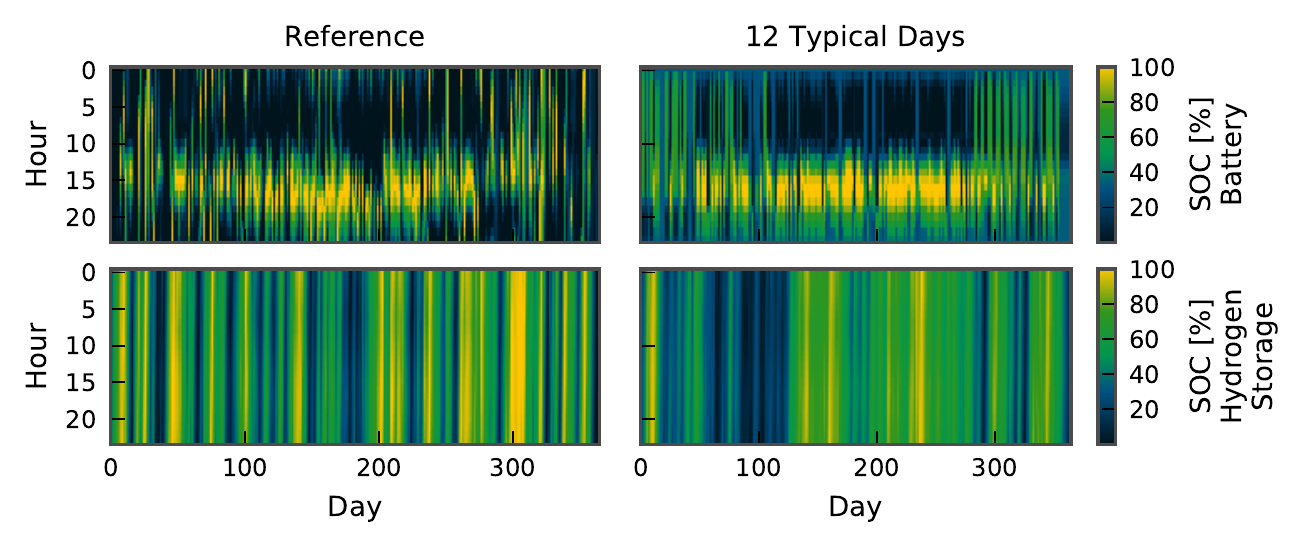}
	\caption{Comparison of the storage operation of the island system based on the full annual time series, in comparison to the storage operation based on 12 typical days and the inter-period state formulation. The colors indicate the state of charge, the x-axis represents the 365 days in the year and the y-axis every 24 hours per day.}
	\label{fig:SOC_TSA}
\end{figure}

In order to validate the new state formulation, the storage inventories of operation with 12 typical periods and the inter-period state formulation is compared to the storage inventory based on the full time series. 

The states of charge of the battery and hydrogen storage are illustrated in Figure \ref{fig:SOC_TSA}. The state of charge for the typical period approach is already presented as additional state of charge of the intra-period and the inter-period states. 

In general, the main operation patterns of the storage system can be captured by the typical period approach with inter-period states. The battery mainly functions as intra-day storage and is getting charged during the day and discharged in the evening. It has more fluctuating operation patterns in the reference case in comparison to the case with 12 typical periods. This is once related to the aggregation which partially reduce the variance of the input profiles, as well due to inter-period states of charge which connect the typical days and limits jumps in the states between the typical days. The hydrogen storage balances weekly wind fluctuations. Nevertheless, because of the reduction of the input data not all operation cycles of the hydrogen storage that were predicted with the full time series can be captured. Still, without the new state description these long term operation cycles could not even have been considered. 

For this example, pressure vessels have been considered as hydrogen storage, and these have a relatively high specific investment. In order to limit their installed capacity, the cost-optimal operation cycles are chosen by the solver for no longer than a few weeks for the reviewed system. In the case large scale hot water storage or caverns as hydrogen storage would be considered, the optimal system operation would shift even more to a long term storage solution with an increased storage capacity. For these types of systems, we expect an even greater necessity to model inter-period states.

\section{Conclusion}
\label{sec:Conclusion}
In the present work, a novel system state description was derived that considers state transitions between typical operation periods which were aggregated with clustering methods. This allows the building of compact temporal energy system models that are still able to consider for the full operational possibilities of long term storage.

The benefit of this method relates to the necessity of long term storage in the considered energy system:
\begin{itemize}
\item For the design of systems where a seasonal storage is no economically feasible option, the new state formulation does not improve the accuracy of the optimization results in comparison to a formulation with independent aggregated periods. Instead, the computational load increases due to the required additional auxiliary variables. Still, a seasonal storage has to be correctly modeled in the solution space to reliably state that it is no economic option, which is not guaranteed with independent periods.
\item	For energy systems which heavily rely on long term storage options, the previous existing approach with independent aggregated typical periods is unable to achieve a system design similar to that based on a full time series operation. The novel introduced state description changes this and a cost optimal system design can be achieved with aggregated typical periods and the information of their sequence, while reducing the overall computational load. 
\end{itemize}

This work resolved a modeling error caused by the formulation of typical periods within energy system models. Nevertheless, the overall results indicate that the solution accuracy is still highly related to the quality of the time series aggregation itself. Therefore, future research should focus on improving the aggregation of design-relevant typical operating periods.

Furthermore, the introduced two-layer state description was used in the example systems to link  typical days within a year, which would also be applicable in cases such as linking typical weeks across a decade. Nevertheless, to also consider energy system design-relevant time series variability within an hour \cite{Stenzel2016} at the same time as variability over the years \cite{Pfenninger2017}, an extension of the two time-layers formulation to a multi-time-layer formulation should be a topic for future analysis.

\section*{Acknowledgments}
\label{s:Acknow}
This work was supported by the Helmholtz Association under the Joint Initiative "EnergySystem~2050 $-$ A Contribution of the Research Field Energy". Further, I want to thank Lara Welder for fruitful discussions about the mathematical derivation of the superposition of states.

\appendix
\label{sec:Appendix}
\section{System modeling}
\label{sec:systemModel}
In order to validate the new model formulation for many different systems with acceptable computational resources, an easily comprehensible simple Mixed-Integer Linear Program has been selected as the system model. More precise models that consider more detailed the size effects of units and investment costs can be found in \citet{Elsido2017, Bahl2017, Gabrielli2017} and \citet{Schutz2017}. 

The system model in this work is defined by a network of specific technologies that are connected by energy flow variables $\dot{E}_{i,j,t}$\nomenclature[A]{$E$}{Energy flow between two components}{}{} at time step $t$. Each connection is therefore defined by an output component $i$ and input component $j$ and belongs to a connection set $(i,j) \in L$. These connections are restricted by the component models  introduced.\nomenclature[B]{L}{Set of component connections}{}{} 

For a typical period, the time steps $t$ are replaced by the time steps $g$ within a single period $k$.

\subsection{Objective function}
The objective function describes the annualized cost of the supply system considered. Therefore, for each device $d$\nomenclature[C]{d}{Considered device or technology}{}{}, the annualized costs are calculated with a capital recovery factor $CRF_{d}$\nomenclature[D]{$CRF$}{Capital Recovery Factor}{}{}, which considers the Weighted Average Cost of Capital \nomenclature[D]{WACC}{Weighted Average Cost of Capital}{}{} $WACC_{d}$ and lifetime $\tau_{d}$\nomenclature[B]{$\tau$}{Lifetime}{}{} of the device in years:
\begin{align}
CRF_{d}=\frac{(1-WACC_{d})^{\tau_{d}} WACC_{d}}{ (1-WACC_{d})^{\tau_{d}}-1}
\end{align}
The capital expenditure of each component is divided into the existing related costs [eur], which only appear if the component is installed, and scale related costs [eur/kW], as well as specific costs, which are scale-dependent \cite{Lindberg2016}. For this reason, each component is modeled by a binary variable $\delta_d$ that defines whether the component exists, and a continuous variable $D_d$, which defines the installed capacity of the component. The resulting device specific annualized fixed cost can be calculated with the existing related capital expenditure ($CAPEX_{exist}$), the scaling-related capital expenditure ($CAPEX_{spec}$)\nomenclature[D]{CAPEX}{Specific capital expenditure}{}{} and fixed operational expenditure ($OPEX_{fix,d}$)\nomenclature[D]{OPEX}{Specific fix operational expenditure}{}{}  as follows 
\begin{align}
\begin{array}{c}
c_{exist,d} =  CAPEX_{exist,d} \left( CRF_{d} + OPEX_{fix,d} \right) \\
c_{spec,d} =  CAPEX_{spec,d} \left( CRF_{d} + OPEX_{fix,d} \right) \\
\end{array}
\end{align}
The costs, which variate with the operation of the system  $c_{var,i,j,t}$, are related to the energy flows $\dot{E}_{i,j,t}$. Along with the scaling of the devices $D_{d}$\nomenclature[A]{$D$}{Scaling of a device}{}{}, the following objective function can be stated:
\begin{align}
\min \sum_{d} c_{exist,d} \delta_{d} + c_{spec,d} D_{d} + \sum_{(i,j) \in L} \sum_{t \in T} c_{var,i,j,t} \dot{E}_{i,j,t} \triangle t
\end{align}

\subsection{Constraints}
The device models establish the constraints of the system. They are divided into five classes, namely: \textit{Source/Sinks, Collectors,  Transformers} and \textit{Storages}. 

The \textit{Source/Sink} class $q$\nomenclature[C]{$q$}{Index of the Source/Sink class}{}{} represents input and output flows to the system, such as photovoltaic feed-in or electricity demand. It is essentially defined by a single equation:
\begin{align}
\eta_{lb,q,t} D_q \leq \sum_{(q,j) \in L } \dot{E}_{q,j,t} \leq \eta_{ub,q,t} D_q \quad \forall \quad t,q
\end{align}
where $\eta_{lb,q,t}D_q$\nomenclature[D]{$LB$}{Lower bound}{}{} \nomenclature[D]{$UB$}{Upper bound}{}{}\nomenclature[B]{$\eta$}{Efficiency}{}{}could be a certain demand that must at least be satisfied at timestep $t$, or $\eta_{ub,q,t}D_q$ could be the maximum photovoltaic feed-in per installed capacity.

The \textit{Collectors} class $n$\nomenclature[C]{$n$}{Index of the Collector class}{}{} can be seen as a hub in which all input energy flows must be equivalent to all output energy flows:
\begin{align}
\sum_{(i,n) \in L } \dot{E}_{i,n,t} - \sum_{(n,j) \in L } \dot{E}_{n,j,t} = 0 \quad \forall \quad t,n
\end{align}

The \textit{Transformer} class $f$\nomenclature[C]{$f$}{Index of the Transformer class}{}{} represents devices that transform the energy from one form to another. Examples include fuel cells or heat pumps. For the definition of these, the energy type (electricity, gas, etc.) $\epsilon$\nomenclature[C]{$\epsilon$}{Energy type}{}{} must be outlined. Each energy flow $\dot{E}_{i,j,t}$ has a certain energy type $\epsilon$. With the energy type's specific transformation efficiency $\eta_{f,\epsilon_{in},\epsilon_{out}}$, the following equation can be stated for each energy transformation in the device:
\begin{align}
\eta_{f,\epsilon_{in},\epsilon_{out}} \sum_{(i,f) \in L,\epsilon_{in} } \dot{E}_{i,f,t} - \sum_{(f,j) \in L,\epsilon_{out} } \dot{E}_{f,j,t} = 0 \quad \forall \quad t,f
\end{align}

The \textit{Storage} class $s$\nomenclature[C]{$s$}{Index of the Storage class}{}{} is defined by an additional variable
the State of Charge $SOC_{s,t}$\nomenclature[A]{$SOC$}{State of charge}{}{} at time step $t$. We can utilize the Euler method to enumerate for the state of charge in the next time step $SOC_{s,t+1}$:
\begin{equation}
\begin{array}{rl}
SOC_{s,t+1} = & SOC_{s,t} (1-\eta^{self}_{s} \Delta t ) \\
& + \eta_{s}^{char} \sum\limits_{ (i,s) \in L} \dot{E}_{i,s,t} \Delta t \\
& -  \frac{1}{\eta^{dis}_{s}} \sum\limits_{ (s,j) \in L } \dot{E}_{s,j,t} \Delta t   \\ 
\end{array}  \quad \forall \quad t,s 
\end{equation}
where $\dot{E}^{char}_{s,t}$ describes the charging flow with an efficiency of $\eta^{char}_{s}$ and $\dot{E}^{dis}_{s,t}$ the discharging flow with related efficiency $\eta^{dis}_{s}$. $\eta^{self}_{s}$ defines the self-discharge of the storage and $\Delta t$\nomenclature[B]{$\Delta t$}{Duration of a single time step}{}{} the step length of a single time step. The state of charge at the beginning of the considered time frame $SOC_{s,1}$ is related to that at the end of the time frame $SOC_{s,N_t+1}$.

The design variable of the storage $s$ which is described by its capacity $D_s$ limits the state of charge to the following:
\begin{equation}
SOC_{s,t} \leq D_s \quad \forall \quad t,s 
\end{equation}

The existing related variable $\delta_d$\nomenclature[A]{$\delta$}{Binary variable determining the existance of a technology}{}{} restricts the scaling-dependent device variable $D_d$ by the so called BigM-Method \cite{Bemporad1999} as follows:
\begin{equation}
\mathbf{M}  \delta_d \geq D_d
\end{equation}
The method is inspired by \citet{Stadler2014} and \citet{Lindberg2016}.

\section{Simplification of the storage operation restrictions}
\label{app:soc_relaxation}
In order to reduce the equation set given by the limitations of the state of charge, shown in equation \ref{eq:SOC_limit}, we introduce two auxiliary variables: $SOC^{intra}_{s,k,max}$ is the maximum state of charge within the typical period $k$ and $SOC^{intra}_{s,k,min}$ is the minimal state of charge. They restrict the intra-period state of charge as follow
\begin{equation}
\begin{array}{c}
SOC^{intra}_{s,k,g} \leq SOC^{intra}_{s,k,max} \quad \forall \quad g,k \\
SOC^{intra}_{s,k,g} \geq SOC^{intra}_{s,k,min} \quad \forall  \quad g,k \\
\end{array}
\end{equation}
The state of charge for the entire sequence of typical periods is then further limited to the maximal and minimal state of charge, as follows:
\begin{equation}
\begin{array}{c}
SOC^{inter}_{s,i} + SOC^{intra}_{s,k=f(i),max}  \leq D_s \quad \forall \quad i \\
SOC^{inter}_{s,i} (1-\eta^{self}_{s} \Delta t )^{N_g} + SOC^{intra}_{s,k,min} \geq 0 \\
 \quad \forall  \quad i 
 \end{array}
\end{equation}
.
The scaling of the storage $D_s$ must to be greater than the high layer state of charge at period $i$ in addition to the maximal state of charge within the related representative period $k$, which can be stated for the equivalent lower bound. 
This modification is a conservative assumption, since the maximum intra-period state of charge is assumed to simultaneously appear to the maximum inter-period state of charge. This holds equivalent for the minimal intra-period inter-period states. Still, for common self-discharge rates, the resulting estimation error engendered by this conservative assumption is negligible.
The number of constraints is thereby reduced from $2xN_ixN_g$ to \mbox{$2x(N_kxN_g+N_i)$} while adding $2xN_g$ additional variables. An application of the approach can be found in \citet{Welder2017}.

\section{Case data}
This section introduces all the relevant data which parameterises the three system models, introduced in Section \ref{sec:Framework}. \ref{sec:TimeSeries} shows the considered raw input time series, while \ref{sec:Parameters} shows the relevant techno-economic input parameters. \ref{sec:scaleResults} shows some example system designs based on typical days and the new state formulation in comparison to the reference system design.

\subsection{Time series data}
\label{sec:TimeSeries}
\begin{figure}[h]
	\centering
  \includegraphics{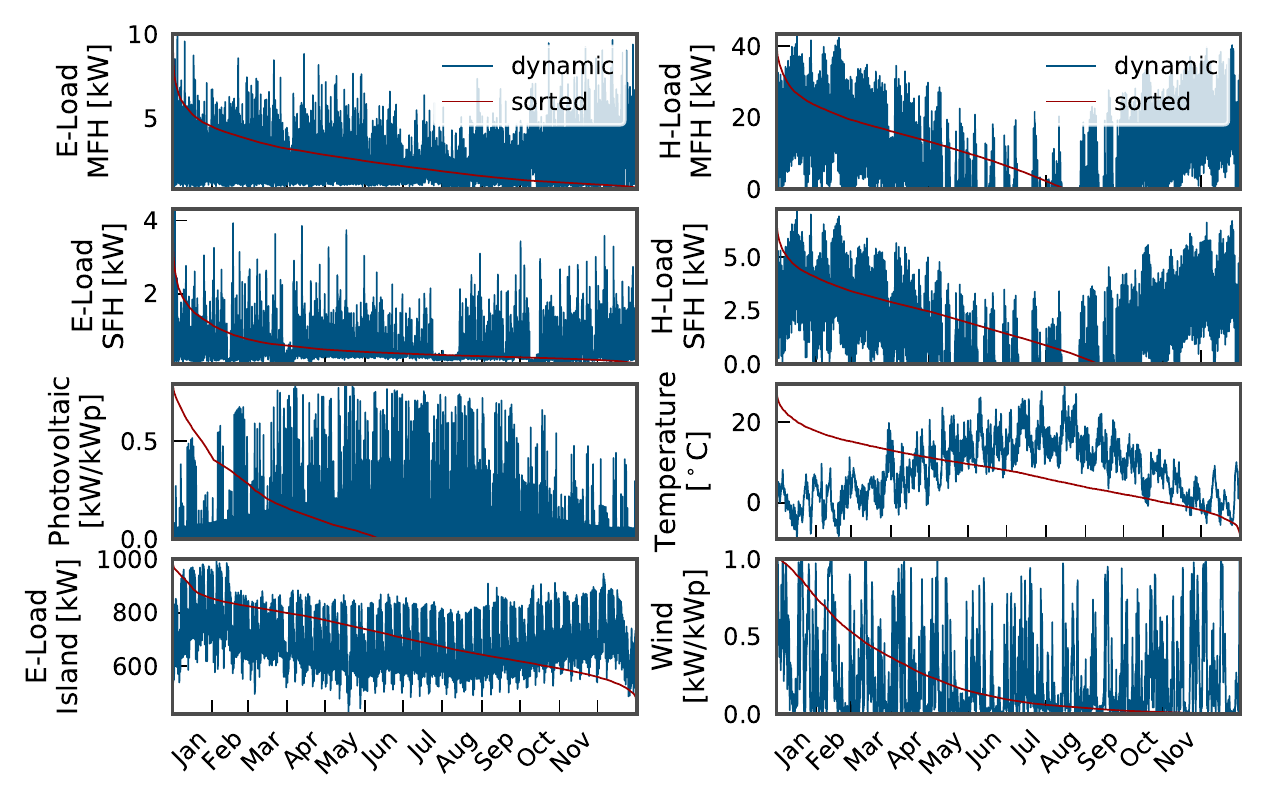}
	\caption{Time series data used for the three case studies as full annual load profiles and as duration curves.}
	\label{fig:Profiles}
\end{figure}

The input time series to the models are shown in Fig. \ref{fig:Profiles}. Their derivation and simulation is explained in Section \ref{sec:Framework}.

\subsection{Technology parameters}
\label{sec:Parameters}

Following technology parameters have been considered for the optimization. The cost approximation is kept simple in order to be able to calculate all cases also with the full time series. The focus of this work was less the research of precise data and more the presentation of practicable methods.  Nevertheless, for the purpose of reproducing the results following data has been used, while the majority of the building technology data is based on \citet{Lauinger2016} and \citet{Lindberg2016a}. The island system is a fictive future case, since the configuration would not be competitive today. The magnitude of cost is based on \cite{DOE2015,DOE2016,IEA2015}.

Table \ref{tab:Converters} includes all data for the considered transformers, the price structure for energy imports or generation units is described in table \ref{tab:Sources}, while table \ref{tab:Storages} includes all storage parameters.

\begin{table}[h]
\caption{Converter parameters}
\label{tab:Converters}
\begin{tabular}{@{}cccccc @{}}
 & \bf{CAPEX} &   & \bf{OPEX} & \bf{Life-} & \bf{Effi- }\\ & exist & spec & fix & \bf{time} & \bf{ciencies}
 \\
\bf{Technology}  & [eur] & [eur/kW$_{el}$] & [\% inv.] & [a]& [-] \\
 \hline
Gas Boiler & 5000 & 50 & 1.5 & 20 & 0.96 \\
CHP & 8000 & 2000 & 5 & 15 & 0.33$\frac{kW_{el}}{kW_{th}}$; 0.52$\frac{kW_{th}}{kW_{th}}$ \\
 \hline
Heat pump & 3000 & 1150 & 2 & 20 & dynamic\cite{Lauinger2016} \\
Electric heater & - &  60 & 0 & 30 & 0.98 \\
\hline
Electrolyser & 100e3 &  500 & 3 & 15 & 0.7 \\
H$_2$ fuel cell & 100e3 &  1100 & 3 & 15 & 0.5 \\
\end{tabular}
\end{table}

\begin{table}[h]
\caption{Source parameters}
\label{tab:Sources}
\begin{tabular}{cccccc}
 & \bf{CAPEX} &   & \bf{OPEX} & &  \bf{Life-} \\
 & exist & spec & fix & var & \bf{time} \\
\bf{Technology}  & [eur] & [eur/kW$_{el}$] & [-] & [eur/kWh] & [a] \\ \hline
Electricity MFH & - & - & 140 eur/a & 0.241 & 50 \\
Gas MFH & - & - & -  & 0.0052 & 50  \\
\hline
Electricity SFH & - & - & -  & 0.281 & 50  \\
Photovoltaic SFH & 1000 & 1200 & 1 \% inv.  & - & 20  \\
\hline
Photovoltaic Island & 10e3 & 800 & 1 \% inv.  & - & 20  \\
Wind turbines & 100e3 & 1000 & 2 \% inv.  & - & 20  \\
Backup plant & - & - & -  & 0.2 & 25  \\
\end{tabular}
\end{table}

\begin{table}[h]
\caption{Storages parameters. The capacity factor of all storage technologies is limited to 1. }
\label{tab:Storages}
\begin{tabular}{@{}ccccccc @{}}
 & \bf{CAPEX}   & \bf{OPEX} & \bf{Life-} &  & \bf{Efficiencies} & \\ &  spec & fix & \bf{time} & charge & discharge & selfdischarge \\
\bf{Technology}  & [eur/kWh] & [\% inv.] & [a]& [-] & [-] &  [1/h]  \\
 \hline
Heat storage & 90 & 0.0 & 25 &  0.99 & 0.99 & 1e-3 \\
\hline
Battery &  300 & 1 & 15 &  0.96 & 0.96 & 5e-4 \\
H$_2$ storage &  15 & 0 & 25 &  0.9 & 1.0 & 0.0 \\
\end{tabular}
\end{table}

\subsection{Selected results}
\label{sec:scaleResults}

\begin{table}
	\centering
\caption{Resulting scaling of the CHP system technologies for different number of typical days, once modeled as independent days and once linked with the new state description, in comparison to the technology scaling based on the full time series.}
\label{tab:scaleCHP}
\begin{tabular}{ll|rrrr}
            &   & CHP &  Gas boiler &  Gas grid &  Heat storage \\
            &  Days &   [kW$_{el}$] &  [kW$_{th}$] &  [kW$_{th}$] & [kWh] \\
\hline
Inde- & 6 &  3.50 &      24.07 &    35.68 &         9.26 \\
   pendent         & 12 &  3.44 &      25.21 &    36.68 &        12.00 \\
            & 21 &  3.39 &      27.24 &    38.66 &        11.82 \\
            & 60 &  3.42 &      27.67 &    39.18 &        17.41 \\
\hline
Linked & 6 &  3.50 &      24.07 &    35.68 &         9.23 \\
            & 12 &  3.44 &      25.32 &    36.79 &        11.76 \\
            & 21 &  3.38 &      27.42 &    38.81 &        11.14 \\
            & 60 &  3.42 &      27.74 &    39.24 &        17.27 \\
\hline
Reference & {} &  3.40 &      27.88 &    39.33 &        16.95 \\
\end{tabular}
\end{table}

\begin{table}
	\centering
\caption{Resulting scaling of the residential system technologies for different number of typical days, once modeled as independent days and once linked with the new state description, in comparison to the technology scaling based on the full time series.}
\label{tab:scaleHP}
\begin{tabular}{ll|rrrr}
\toprule
            &   &  ElectricHeater &  HeatPump &  HeatStorage &  Photovoltaic \\                   
            &Days  &          [kW] &     [kW] &             [kWh] &        [kW] \\ 
\hline
Inde- & 6 &            0.00 &      3.89 &        10.40 &          2.27 \\
    pendent        & 12 &            0.12 &      3.81 &        10.18 &          3.43 \\
            & 21 &            0.45 &      3.69 &         9.87 &          3.36 \\
            & 60 &            0.46 &      3.77 &        10.09 &          3.39 \\ \hline
Linked & 6 &            0.00 &      3.89 &        10.35 &          2.27 \\
            & 12 &            0.00 &      3.81 &        10.14 &          3.18 \\
            & 21 &            0.35 &      3.79 &         9.79 &          3.30 \\
            & 60 &            0.48 &      3.77 &        10.10 &          3.37 \\ \hline
Reference & {} &            0.67 &      3.76 &        12.38 &          3.43 \\
\end{tabular}
\end{table}

\begin{table}
	\centering
\caption{Resulting scaling of the Island system technologies for different number of typical days, once modeled as independent days and once linked with the new state description, in comparison to the technology scaling based on the full time series.}
\label{tab:scaleIsland}
\begin{tabular}{ll|rrrrrrr}
        &    &  \makecell{Backup \\ plant} &  Battery &    \makecell{Electro- \\ lyser}&  \makecell{Fuel \\ cell} & \makecell{H$_2$ \\ storage} &  \makecell{Photo \\ voltaic} & \makecell{ Wind \\ turbine} \\
                   &Days  &          [kW] &     [kWh] &             [kW] &        [kW] &               [kWh] &          [kW] &  [kW] \\    
\hline
Inde- & 6 &          231 &     4863 &             0 &         0 &                0 &          5366 &  1726 \\
   pendent         & 12 &          501 &     7123 &             0 &         0 &                0 &          5335 &  1923 \\
            & 21 &          486 &     6726 &             0 &         0 &                0 &          4879 &  2489 \\
            & 60 &          613 &     6646 &             0 &         0 &                0 &          4781 &  2689 \\
\hline
Linked & 6 &          232 &     5372 &             0 &         0 &                0 &          4965 &  1952 \\
            & 12 &          385 &     3502 &           654 &       207 &            51917 &          3694 &  2715 \\
            & 21 &          408 &     3719 &           625 &       187 &            71115 &          3534 &  2713 \\
            & 60 &          428 &     3104 &           804 &       256 &            62309 &          3236 &  2783 \\
\hline
Reference & {} &          377 &     3364 &           749 &       256 &            42663 &          3329 &  2633 \\
\end{tabular}
\end{table}


\clearpage

\section*{References}

\bibliography{mybibfile}

\end{document}